%% file: main.tex
\documentclass[preprint]{jpsj3}
\usepackage{txfonts}

\usepackage[whole]{bxcjkjatype}  

\usepackage{xcolor}
\usepackage{url}
\usepackage{ulem}

\input{math_commands.tex}

\DeclareRobustCommand{\remark}[2]{\underset{#1}{\uwave{#2}}}

\DeclareRobustCommand{\a}[0]{{\bm{a}}}
\DeclareRobustCommand{\b}[0]{{\bm{b}}}

\DeclareRobustCommand{\n}[0]{{\bm{n}}}

\DeclareRobustCommand{\w}[0]{{\bm{w}}}
\DeclareRobustCommand{\bwt}[0]{{\tilde{\bm{w}}}}
\DeclareRobustCommand{\wt}[0]{\tilde{w}}

\DeclareRobustCommand{\wsfss}[0]{\hat{\mathsf{w}}_{\rm SS}}
\DeclareRobustCommand{\wsfko}[0]{\hat{\mathsf{w}}_{\rm dKO}}
\DeclareRobustCommand{\wtsfko}[0]{\hat{\tilde{\mathsf{w}}}_{\rm dKO}}

\DeclareRobustCommand{\hatqss}[0]{\hat{Q}_{\rm SS}}
\DeclareRobustCommand{\hatmss}[0]{\hat{m}_{\rm SS}}
\DeclareRobustCommand{\hatchiss}[0]{\hat{\chi}_{\rm SS}}
\DeclareRobustCommand{\hatvss}[0]{\hat{v}_{\rm SS}}

\DeclareRobustCommand{\hatqko}[0]{\hat{Q}_{\rm dKO}}
\DeclareRobustCommand{\hatqkko}[0]{\hat{\tilde{Q}}_{{\rm dKO}}}
\DeclareRobustCommand{\hatmko}[0]{\hat{m}_{\rm dKO}}
\DeclareRobustCommand{\hatchiko}[0]{\hat{\chi}_{\rm dKO}}
\DeclareRobustCommand{\hatvko}[0]{\hat{v}_{\rm dKO}}
\DeclareRobustCommand{\hatvkko}[0]{\hat{\tilde{v}}_{{\rm dKO}}}

\DeclareRobustCommand{\dalgo}[0]{D^\ast}

\DeclareRobustCommand{\p}[0]{\partial}

\title{Replica Analysis for Ensemble Techniques in Variable Selection}

\author{Takashi Takahashi\thanks{takashi-takahashi@g.ecc.u-tokyo.ac.jp}}
\inst{Institute for Physics of Intelligence, the University of Tokyo, 7-3-1 Hongo, Bunkyo-ku, Tokyo, 113-0033, Japan} 

\abst{
    Variable selection is a problem of statistics that aims to find the subset of the $N$-dimensional possible explanatory variables that are truly related to the generation process of the response variable. In high-dimensional setups, where the input dimension $N$ is comparable to the data size $M$, it is difficult to use classic methods based on $p$-values. Therefore, methods based on the ensemble learning are often used. In this review article, we introduce how the performance of these ensemble-based methods can be systematically analyzed using the replica method from statistical mechanics when $N$ and $M$ diverge at the same rate as $N,M\to\infty, M/N\to\alpha\in(0,\infty)$. As a concrete application, we analyze the power of stability selection (SS) and the derandomized knockoff (dKO) with the $\ell_1$-regularized statistics in the high-dimensional linear model. The result indicates that dKO provably outperforms the vanilla knockoff and the standard SS, while increasing the bootstrap resampling rate in SS might further improve the detection power.
}

\begin{document}
\maketitle

\section{Introduction}
\label{section: introduction}

Variable selection refers to a statistical task of estimating which of the many explanatory variables, potentially related to the data generation process, are truly related to it. Such tasks are often encountered in various scientific domains ranging from astronomy to economics \cite{uemura2015variable, uemura2016data, Luca2015big, igarashi2018exhaustive}. However, because the data is usually of limited size and contains some noise, it is difficult to identify the truly relevant explanatory variables without an error, except in some special situations \cite{wainwright2009sharp,okajima2023average}. Therefore, it is important to quantify the importance of each explanatory variable and to perform variable selection while controlling the false discovery rate (FDR), where FDR is defined as the rate of falsely identified variables among those estimated to be related.

The classic variable selection method that can control FDR is the hypothesis testing based on $p$-values. Once $p$-values are evaluated for each explanatory variables, important variables can be determined while  controlling the FDR using testing methods such as the Benjamini-Hochberg algorithm \cite{benjamini1995controlling}. However, in modern statistical problems, the number of the potential explanatory variables is often comparable to the data size. In such high-dimensional setups, evaluation of $p$-values is non-trivial and no versatile algorithm that can yield exact $p$-values is known although some approaches are proposed \cite{zhang2013confidence,geer2014asymptotically, javanmard2014confidence,javanmard2014hypothesis, taylor2015statistical,Takahashi_2018,sawaya2024high}.

Over the last decade or so, instead of using such a classic hypothesis testing, ensemble learning-like approaches have been proposed \cite{meinshausen2010stability,shah2012variable,tuan2020aggregation,Ren2023derandomizing,ren2023derandomised,dai2023false,dai2023scale}. The basic idea of this approach is to repeatedly apply some variable selection algorithm, which itself may not guarantee control the FDR, to data with additional randomness added manually, calculate the selection probability for each variable, and then obtain the final result based on those probabilities.  Methods of this type include stability selection (SS) using resampled data \cite{meinshausen2010stability,shah2012variable} and the derandomized knockoff (dKO) \cite{tuan2020aggregation, ren2023derandomised,Ren2023derandomizing}, which will be the focus of analysis in this paper. See Sec. \ref{section: setup} for details of these algorithms. These methods are easily applicable even in high-dimensional settings, can control FDR, despite the selection probability is not equal to traditional $p$-values.

Although both SS and dKO (and their variants) can control FDR, they actually sacrifice the effective sample size in different ways to assess selection probability. Thus, quantitative performance evaluation of these methods is important to provide guidance on which method to use in a given situation. Performance measures of interest include power,  defined as the true positive rate (TPR) at a given FDR, where TPR is defined as the rate of truly identified variables. However, in many cases, such a performance comparison is limited to experimental evaluations on synthetic data, and the theoretical evaluation is scarce.

Given the above situation, this paper presents a method for systematically characterizing the performance of these ensemble learning-based variable selection methods using the technique of statistical physics and aims to compare the performance of SS and dKO in a simplified setting, in the limit where the data size $M$ and the number of potential explanatory variables $N$ diverges at the same rate as $N,M\to\infty, M/N\to\alpha\in(0,\infty)$. In the following, $N\to\infty$ is used to represent this proportional limit as a shorthand notation.

There are two salient features in this asymptotic regime. First, the macroscopic quantities such as TPR and FDR do not depend on the details of the realization of the training data in this asymptotic regime. That is, the fluctuations of these quantities with respect to the training data vanish, allowing us to make sharp theoretical predictions. This is analogous to macroscopic physical systems where we can make sharp predictions about the macroscopic quantities such as specific heat, pressure, etc., and this proportional asymptotic regime is sometimes referred to as the thermodynamic limit (TDL), following the custom of physics even in the context of statistics. Second, the performance of the algorithm can be fully characterized by a few macroscopic quantities if the data generation process and the algorithm are simple enough. This is due to the fact that the statistical problems are solvable in simple settings in a manner similar to infinite-range Ising models in statistical physics \cite{nishimori2001statistical, nishimori_elements_2010, charbonneau2023spin}, and it suggests that mean-field approaches can also be useful as a starting point for an analysis even in statistics.

The analysis of statistics and machine learning methods in TDL has been actively carried out since the 1990s, using techniques from statistical physics \cite{engel2001statistical,lenka2016statistical,montanari2018mean,miolane2019fundamental, barbier2020high, barbier2020mean, sakata2023prediction, montanari2024friendly}. Recently, it has become one of the standard analytical approaches, and ensemble learning methods have also been studied in the same proportional asymptotics. The literature \cite{Sollich1995learning, Krogh1997statistical, malzahn2001variational, malzahn2002statistical, malzahn2003learning, lejeune2020implicit, dascoli2020double, loureiro2022fluctuations, Takahashi2023role, clarte2024analysis, patil2024asymptotically, takahashi2024a}
analyze the performance of the ensemble average of linear models. In addition, the literature \cite{obuchi2019semi, takahashi2019replicated, Takahashi2020semi} develops approximate algorithm for bootstrap methods \cite{efron1992bootstrap, bach2008bolasso,meinshausen2010stability} that can avoid repeated fitting to resampled data by combining the replica method of statistical physics and approximate inference algorithms of machine learning. However, the focus of these studies has been limited to evaluating generalization performance in the context of machine learning or reducing the computational burden of algorithms, and performance evaluation in the context of mathematical statistics has rarely been performed. 

In this paper, we present the approach based on the replica method, which is the method of statistical physics developed to analyze the physics of disordered systems such as the spin glasses \cite{mezard1987spin,fischer1993spin,charbonneau2023spin} and structural glasses \cite{binder2011glassy,parisi2020theory}. We will explain how the ensemble learning-based methods can be analyzed using the replica method in a rather general framework, and then apply it to the analysis of SS and dKO with $\ell_1$ regularized statistics.

The remainder of the paper is organized as follows. Sec. \ref{section: setup} explains the problem formulation. In Sec. \ref{section: replica method}, we describe the replica method for analyzing ensemble learning-based methods, where the statistical mechanics formulation of the inference problem is also introduced. Then, in Sec. \ref{section: mean-field picture}, we present the mean-field picture of SS and dKO, where their behavior is described by an effective single-body problem with parameters determined as solutions to self-consistent equations. Based on this mean-field picture, the performance analysis of SS and dKO is presented in Sec. \ref{section: performance comparison}. Finally, Sec. \ref{section: summary and discussion} provides a summary with some concluding discussions.

\begin{table}[t!]
    \centering
    \begin{tabular}{|l l|}
    \hline
    Notation & Description \\
    \hline\hline
        $[n]$ & for a positive integer $n$, the set $\{1,\dots,n\}$
        \\
        $\bm{x} \cdot \bm{y}$ & for vectors $\bm{x},\bm{y}\in\R^N$, the inner product of them: 
        \\
        & $\bm{x} \cdot \bm{y}=\sum_{i=1}^N x_iy_i$.
        \\
        $\1(\cdot)$    & indicator function
        \\
        $\gN(\mu, \sigma^2)$ & Gaussian density with mean $\mu$ and variance $\sigma^2$
        \\
        ${\rm Poi}(\mu_\ast)$ & Poisson distribution with mean $\mu_\ast$
        \\
        iid & independent and identically distributed
        \\
        $\E_{X\sim p_X}[f(X)]$  & Expectation regarding random variable $X$ where $p_X$ is
        \\
        &  the density function for the random variable $X$
        \\
        & (lower subscript $X \sim p_X$ or $p_X$ can be omitted 
        \\
        & if there is no risk of confusion)
        \\
        $|S|$ & for a set $S$, the size of a set $S$.
        \\
        \hline
    \end{tabular}
    \caption{Notations}
    \label{table: notations}
\end{table}

\subsection{Notations}
\label{subsec: notations}

Throughout the paper, we use some shorthand notations for convenience. We summarize them in Table \ref{table: notations}. 

\section{Setup}
\label{section: setup}
This section presents the problem formulation. The assumptions about the data generation process are first described, and then the SS and dKO procedures are formalized. 

\subsection{Data generation process}
Let $D=\{(\bm{x}_\mu, y_\mu)\}_{\mu=1}^M, \bm{x}_\mu\in\R^N, y_\mu \in \R$ be the set of independent and identically distributed (iid) data points, where $\bm{x}_\mu$ be the vector of potential explanatory variables and $y_\mu$ be the corresponding response. We consider the simplest model for the generation of the response $y_\mu$, i.e., the noise-corrupted linear measurement:
\begin{equation}
    y_\mu = \bm{x}_\mu \cdot \bm{w}_0 + \epsilon_\mu, \quad \epsilon_\mu \sim \gN(0,\Delta),
\end{equation}
where $\epsilon_\mu$ is the measurement noise, and $\bm{w}_0\in\R^N$ is the true regression coefficient. We also assume that $\bm{x}_\mu \sim_{\rm iid} \gN(0,I_N/N)$ as the simplest model of the random generation of the explanatory variable.

The variable selection problem is motivated by the belief that, in many practical applications, the response $y_\mu$ depends on only a subset of $\bm{x}_\mu$. To implement this idea, we assume that $\bm{w}_0$ is sparse in the sense that the number of non-zero elements of $\bm{w}_0$ is limited to $\rho N, ~ \rho\in [0,1]$. Specifically, we assume that each component of $\bm{w}_0$ is generated from the Gauss-Bernoulli distribution:
\begin{equation}
    w_{0,i} \sim_{\rm iid} p_{w_0}, \quad p_{w_0}(\cdot) = \rho\phi_+ + (1-\rho)\delta_0,
\end{equation}
with $\phi_+=\gN(0,1)$.  Then, the goal of the variable selection is to find the support $\gS=\{i \mid  w_{0,i}\neq 0, i = 1,\dots,N\}\subset \{1,\dots,N\}$ from $D$ since $w_{0,i}\notin\gS$ clearly does not contribute to the generation of $y_\mu$. 

\subsection{Algorithms}
To estimate the support $\gS$, we consider SS and dKO that are described below. Both algorithms use a randomized variable selection algorithm to compute the selection probability $\Pi_i(D)$ for each explanatory variable, with respect to the randomness induced in the randomized algorithm, and then, give the estimator of $\gS$ as the set of indices that has higher selection probability than some threshold $\Pi_{\rm th}$ as 
\begin{equation}
    \hat{\gS}(D) = \{i \mid \Pi_i(D)>\Pi_{\rm th}, i\in[N]\},
\end{equation}
where $[N]=\{1,2,\dots,N\}$. Although they can be considered as similar methods, they use rather different randomized algorithms to calculate the selection probability.

SS uses $\ell_1$ regularized linear regression, also known as the lasso \cite{tibshirani1996regression}, for bootstrap data. Specifically, bootstrap data $D^\ast_{\rm SS}=\{(\bm{x}_\mu^\ast, y_\mu^\ast )\}_{\mu=1}^{M_B}, (\bm{x}_\mu^\ast, y_\mu^\ast )\in D$ is generated by sampling $M_B$ point from $D$ with replacement, and then apply $\ell_1$ regularized linear regression to $D_{\rm SS}^\ast$ as 
\begin{align}
    \hat{\bm{w}}_{\rm SS}(D,D_{\rm SS}^\ast) &=\argmin_{\bm{w}\in\R^N} \gL_{\rm SS}(\bm{w};D,D_{\rm SS}^\ast),
    \label{eq: estimator SS}
    \\
    \gL_{\rm SS}(\bm{w};D,D_{\rm SS}^\ast) &= \sum_{\mu=1}^M \frac{1}{2}(y_\mu^\ast - \bm{x}_\mu^\ast \cdot \bm{w})^2 + \lambda\sum_{i=1}^N |w_i|
    \label{eq: cost function SS}
\end{align}
where $\lambda>0$ is the regularization parameter. Since the $\ell_1$ regularization can select the variables by making some component of $\hat{\bm{w}}_{\rm SS}$ exactly zero, we can compute the selection probability as 
\begin{equation}
    \Pi_{{\rm SS},i}(D) = \E_{D_{\rm SS}^\ast}\left[
        \1\left(
            \hat{w}_{{\rm SS},i}(D,D_{\rm SS}^\ast)\neq 0
        \right)
    \right].
\end{equation}

dKO uses the knockoff (KO) algorithm \cite{barber2015controlling, candes2018panning} as its randomized variable selection algorithm. Specifically, in KO, we randomly generate fake explanatory variable $\tilde{\bm{x}}_\mu, \mu\in[M]$ that is statistically similar to $\bm{x}_\mu$. In our setup, generating $\tilde{\bm{x}}_\mu\sim\gN(0,I_N/N)$ is sufficient. Then, using both $\tilde{\bm{x}}_\mu$ and $\bm{x}_\mu$, we perform $\ell_1$ regularized linear regression as 
\begin{align}
    &\left(
        \hat{\bm{w}}_{\rm dKO}, \bwt_{\rm dKO}
    \right)(D,\tilde{X})
    = \argmin_{(\bm{w}, \bwt)\in\R^{2N}} \gL_{\rm dKO}((\bm{w}, \bwt);D,\tilde{X})
    \label{eq: estimator dKO}
    \\ 
    &\begin{aligned}
    \gL_{\rm dKO}((\bm{w}, \bwt);D,\tilde{X}) = \sum_{\mu=1}^M &\frac{1}{2}(y_\mu - \bm{x}_\mu\cdot\w - \tilde{\bm{x}}_\mu\cdot\bwt)^2 ,
        \\ 
        &+ \lambda\sum_{i=1}^N (|w_i| + |\wt_i|)
    \end{aligned}
    \label{eq: cost function dKO}
\end{align}
where $\tilde{X} = [\tilde{x}_{\mu i}]\in\R^{M\times N}$. The idea of KO is to select variable $i$ if $Z_i(\tilde{X}) = |\hat{w}_{{\rm dKO}, i}| - |\hat{\wt}_{{\rm dKO},i}|$  is larger than a threshold $Z_{\rm th}$, where $Z_i$ is called the lasso coefficient difference (LCD) statistics, based on the intuition that $\tilde{x}$ is irrelevant to the generation of the response $y_\mu$. Since LCD fluctuates with respect to $\tilde{X}$, dKO computes the selection probability as 
\begin{equation}
    \Pi_{{\rm dKO}, i}(D) = \E_{\tilde{X}}\left[
        \1\left(Z_i(\tilde{X}) > Z_{\rm th}\right)
    \right].
\end{equation}

As mentioned in Sec. \ref{section: introduction}, SS and dKO sacrifice the effective sample size in different ways to assess the selection probability. In SS, because $D_{\rm SS}^\ast$ contains multiple instances of the same data point, the number of unique data point in $D_{\rm SS}^\ast$ is smaller than $M$ as long as $M_B/M<\infty$, which effectively reduces the sample size of each bootstrap data. In dKO, the number of the regression coefficient is doubled because the knockoff variable $\tilde{\bm{x}}_\mu\in\R^N$ is introduced, which effectively reduces the sample size compared to the number of the model parameter. There is an intuitive argument \cite{Ren2023derandomizing} that increase in parameter dimension is better than directly decreasing unique data points. However, it does not explain how much performance difference would be observed quantitatively until the values are evaluated precisely. Showing this difference qualitatively is one goal of this paper.

Specifically, we will mainly focus on the evaluation of the power, i.e., TPR at a given FDR. In this setup, TPR and FDR are given as follows:
\begin{align}
    {\rm TPR} &= \frac{
        \left|
            \{i \mid i\in [N], \Pi_i > \Pi_{\rm th}, w_{0,i}\neq 0\}
        \right|
    }{
        |\{i \mid i\in[N], w_{0,i}\neq 0\}|
    },
    \\
    {\rm FDR} &= \frac{
        \left|
            \{i \mid i\in [N], \Pi_i > \Pi_{\rm th}, w_{0,i} =  0\}
        \right|
    }{
        \left| \hat{\gS}(D)\right|
    }.
\end{align}

\section{Replica Method}
\label{section: replica method}

\begin{figure}
    \centering
    \includegraphics[width=1.0\linewidth]{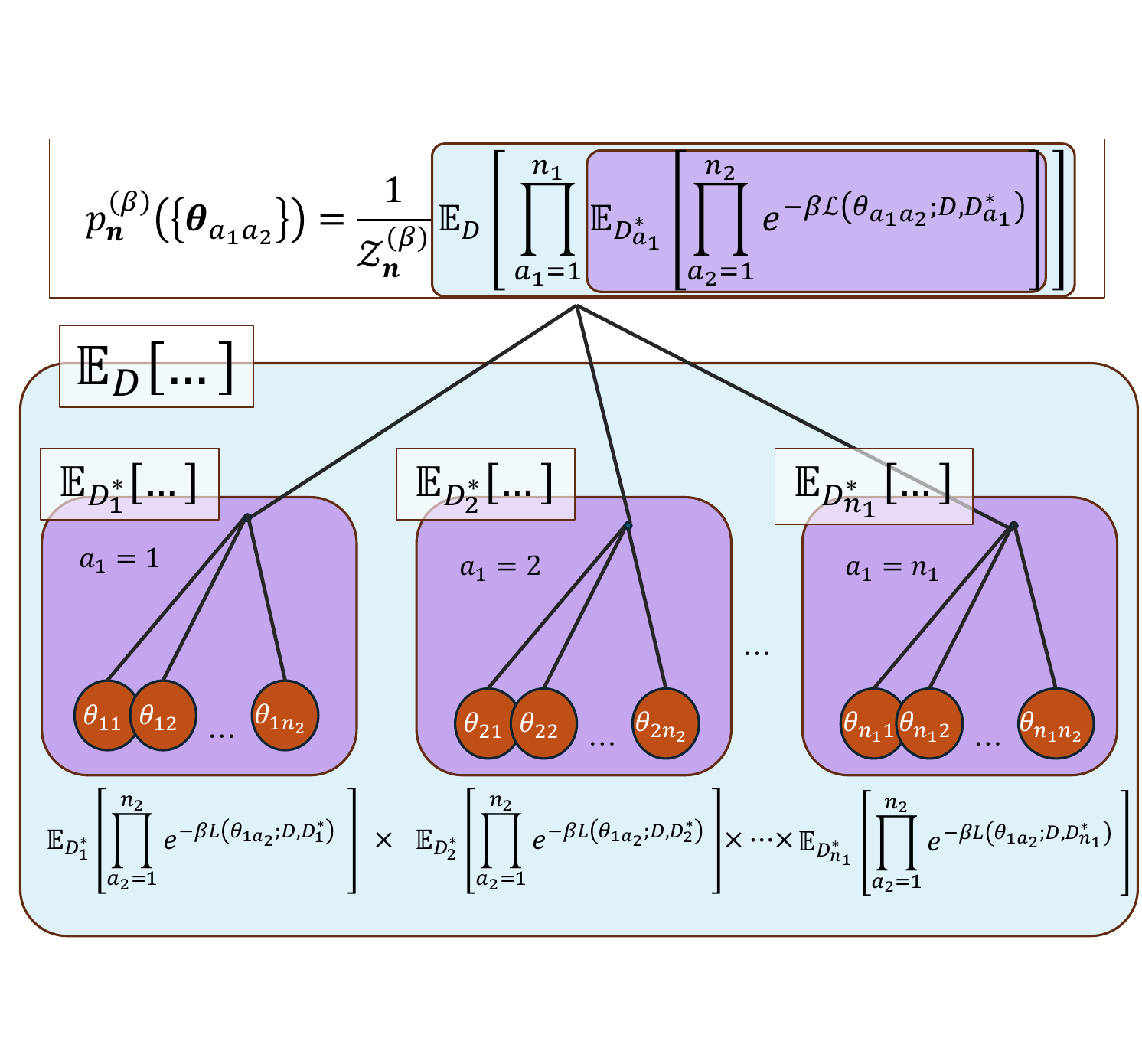}
    \caption{
        Graphical representation of the replicated system defined as \eqref{eq: replicated system}.
    }
    \label{fig: replicated system}
\end{figure}

To evaluate the performance of SS and dKO, we need to characterize how each component of the estimators \eqref{eq: estimator SS} and \eqref{eq: estimator dKO} fluctuate statistically over the randomness in the acquired data $D$ or that introduced by the algorithms $D_{\rm SS}^\ast, \tilde{X}$. In this section, we introduce an approach to this problem based on the replica method from statistical mechanics.

In the following, to describe SS and dKO in a unified manner, let $\bm{\theta}$ be the parameters to be estimated, and let $\dalgo$ be the randomness introduced by the algorithm. That is, $\bm{\theta}, \dalgo$ represents $\bm{w}\in\R^N, D_{\rm SS}^\ast$ in SS, and $(\bm{w}, \bwt)\in\R^{2N}, \tilde{X}$ in dKO, respectively. Furthermore, let $\gL(\bm{\theta};D,\dalgo)$ be the cost function to be optimized, i.e., $\gL(\bm{\theta};D,\dalgo)$ represents $\gL_{\rm SS}$ in \eqref{eq: cost function SS} in SS, and $\gL_{\rm dKO}$ in \eqref{eq: cost function dKO}, respectively.

\subsection{Statistical mechanics formulation}
The first step is to formulate the optimization of random cost functions \eqref{eq: estimator SS} and \eqref{eq: estimator dKO} as a problem of statistical mechanics of disordered system. 

For this, we introduce a probability density $p_{\rm B}^{(\beta)}(\bm{\theta} \mid D,\dalgo)$, which is defined as 
\begin{equation}
    p_{\rm B}^{(\beta)}(\bm{\theta}\mid D,\dalgo) = \frac{1}{Z^{(\beta)}(D,\dalgo)}e^{-\beta \gL(\bm{\theta};D,\dalgo)},
    \label{eq: boltzmann dist general}
\end{equation}
where $\beta\in(0,\infty)$ is the positive parameter termed the inverse temperature following the custom of statistical physics, and $Z^{(\beta)}(D,\dalgo) = \int e^{-\beta \gL(\bm{\theta};D,\dalgo)}d\bm{\theta}$ is the normalization constant termed the partition function.  Here, if the cost function is regarded as a random Hamiltonian with a dynamical variable $\bm{\theta}$ and a quenched disorder $D,D^\ast$, then, $p_{\rm B}^{(\beta)}$ can be seen as the Boltzmann distribution at an inverse temperature $\beta$. At $\beta\to\infty$, the Boltzmann distribution \eqref{eq: boltzmann dist general} converges to the uniform distribution over the minimum of $\gL(\bm{\theta};D,\dalgo)$. Since the estimators \eqref{eq: estimator SS} and \eqref{eq: estimator dKO} are given as the minimizer of $\gL(\bm{\theta};D,\dalgo)$, analyzing these estimators is equivalent to analyzing the zero-temperature limit of the Boltzmann distribution.

The above statistical mechanics formulation is a formal rewrite and does not, in itself, reduce the difficulty of the problem. However, rewriting in this way naturally allows us to introduce the techniques from statistical mechanics, such as the replica method.

\subsection{Replicated system}

The Boltzmann distribution \eqref{eq: boltzmann dist general} is for describing how $\bm{\theta}$ fluctuates \textit{thermally} for a fixed quenched randomness $D, D^\ast$ at a given inverse temperature $\beta$. Therefore, to consider the fluctuation of the ground state in $D,D^\ast$, a special method must be used for this purpose. In particular, since the algorithmic randomness $D^\ast$ is generated conditioned by the acquired data $D$, there is a hierarchical structure of the randomness. Therefore, we need to characterize the fluctuations of $D$ and $D^\ast$, separately.

The method for handling such kind of the fluctuation is the replica method \cite{mezard1987spin, charbonneau2023spin}. The replica method uses the probability density, which we refer to as the \textit{replicated system}.  Let $n_1,n_2\in\sN$ be the positive integers. Then, the replicated system is defined as follows
\begin{equation}
    p^{(\beta)}_{\n}(\{\bm{\theta}_{\a}\}) = \frac{1}{\gZ^{(\beta)}_{\n}}\E_D\left[
        \prod_{a_1=1}^{n_1}\E_{D^\ast_{a_1}}\left[
            \prod_{a_2=1}^{n_2}e^{-\beta \gL(\bm{\theta}_{\a};D,D^\ast_{a_1})}
        \right]
    \right],
    \label{eq: replicated system}
\end{equation}
where $\a=(a_1,a_2), a_1\in[n_1],a_2\in[n_2]$, which is termed the replica indices, is the shorthand notation for the subscript, and $\gZ_{\n}^{(\beta)}=\E_D[\E_{D^\ast}[Z(D,D^\ast)^{n_2}]^{n_1}]$ is the normalization constant.  This is a density over $\R^{n_1\times n_2 \times {\rm dim}\bm{\theta}}$. 

Note that the replicated system is not a conditional distribution of the quenched randomness anymore because the average over $D,\{D^\ast_{a_1}\}$ is taken explicitly. Due to this averaging, $\bm{\theta}_\a$ with different replica indices are not independent. For example, 
\begin{equation}
    \langle
        \theta_{\a,i}\theta_{\b,i}
    \rangle_{\n}\neq 
    \left\langle \theta_{\a,i}\right\rangle_\n
    \left\langle \theta_{\b,i}\right\rangle_\n, 
    \quad \a\neq\b,
\end{equation}
where $i$ is an index to distinguish the component of the parameter $\bm{\theta}$ and $\langle \dots \rangle_\n$ represents the average over the replicated system \eqref{eq: replicated system}. Furthermore, since $D_{a_1}^\ast$ has a replica index, the replicated system has a trivial hierarchical structure (see Fig. \ref{fig: replicated system} for a graphical representation on this point), which is crucial for characterizing the fluctuations of $D$ and $D^\ast$ separately, as we will see below.

The idea of the replica method is to obtain some useful information regarding the fluctuation over $D,D^\ast$ by formally continuing some nontrivial analytical results on the replicated system defined on positive integers $n_1,n_2\in\sN$ to $n_1, n_2\in\R$. To get an intuition about this, let us see the $K$-point correlation function 
$
\langle
    \theta_{\a^{(1)},i} \theta_{\a^{(2)},i}  \dots \theta_{\a^{(K)},i} 
\rangle_\n,
~
K\in\sN
$.
By a formal computing, this can be written, depending on the choice of the replica index $\{\a^{(k)}\}_{k=1}^K$, as follows:
\begin{align}
    &\langle
        \theta_{\a^{(1)},i} \theta_{\a^{(2)},i}  \dots \theta_{\a^{(k)},i} 
    \rangle_\n
    \nonumber 
    \\
    &=\left\{
            \begin{array}{l}
                \E_D\left[
                    \frac{
                        \E_{D^\ast}[\langle \theta_i\rangle^K Z^{(\beta)}(D,D^\ast)^{n_2}]
                    }{
                        \E_{D^\ast}[Z^{(\beta)}(D,D^\ast)^{n_2}]
                    }
                    \E_{D^\ast}[Z^{(\beta)}(D,D^\ast)^{n_2}]^{n_1}
                \right],
                \vspace{1truemm}
                \\
                \hfill
                a_1^{(1)} =\dots=a_1^{(K)}, a_2^{(1)}\neq\dots,\neq a_2^{(K)}
                \vspace{2truemm}
                \\ 
                \E_D\left[
                    \left(
                        \frac{
                            \E_{D^\ast}[\langle \theta_i\rangle Z^{(\beta)}(D,D^\ast)^{n_2}]
                        }{
                            \E_{D^\ast}[Z^{(\beta)}(D,D^\ast)^{n_2}]
                        }
                    \right)^K
                    \E_{D^\ast}[Z^{(\beta)}(D,D^\ast)^{n_2}]^{n_1}
                \right],
                \vspace{1truemm}
                \\
                \hfill
                a_1^{(1)}\neq \dots \neq a_1^{(K)},
            \end{array}
    \right.
    \label{eq: K-point correlation at finite n}
\end{align}
where $\langle \dots \rangle$ (bracket without lower subscript $\n$) is an expectation over the Boltzmann distribution \eqref{eq: boltzmann dist general}. Since the above expression depends on $n_1, n_2$ in the form of the power function, we can consider formal analytical continuation of $n_1,n_2\in\sN$ to $n_1, n_2\in\R$, and take the limit $n_1,n_2\to0$. Then, we obtain the following expression:
\begin{align}
    &\langle
        \theta_{\a^{(1)},i} \theta_{\a^{(2)},i}  \dots \theta_{\a^{(k)},i} 
    \rangle_\n
    \nonumber 
    \\
    &\overset{n_1,n_2\to0}{\longrightarrow} 
    \left\{
        \begin{array}{l}
            \E_{D,D^\ast}\left[
                \langle \theta_i\rangle^K
            \right],
            \vspace{1truemm}
            \\
            \hfill
            a_1^{(1)} =\dots=a_1^{(K)}, a_2^{(1)}\neq\dots,\neq a_2^{(K)}
            \vspace{2truemm}
            \\ 
            \E_{D}\left[
                \E_{D^\ast}\left[
                    \langle \theta_i\rangle
                \right]^K
            \right],
            \vspace{1truemm}
            \\
            \hfill
            a_1^{(1)}\neq \dots \neq a_1^{(K)}.
        \end{array}
    \right.
\end{align}
This implies that, in the replicated system \eqref{eq: replicated system}, the correlation between elements of $\{\bm{\theta}_\a\}$ with different first replica indices $a_1$, i.e., replicas chosen from different purple boxes in Fig. \ref{fig: replicated system}, encodes the fluctuation of $D^\ast$ conditioned by $D$. Similarly, the correlation between points with different second replica indices $a_2$, i.e., replicas chosen from a same purple box in Fig. \ref{fig: replicated system}, encodes the fluctuation of both $D^\ast, D$. 

Although the above computations are completely formal and require solving a multi-body problem, the point is that, when using a simple data model, the behavior of the replicated system can be fully characterized by an effective single-body problem. This single-body problem is characterized by a small finite number of order parameters, and the randomness arising from $D,D^\ast$ are also effectively characterized by one-dimensional random variables, which drastically simplifies the problem in high dimension. Furthermore, the analysis of the replicated system is completely analogous to the one-step replica symmetry breaking (1RSB) solution of the spin-glass systems, reflecting the trivial hierarchical structure of the replicated system. Thus, no novel idea is required for the analysis of the replicated system itself, as long as no additional replica symmetry breaking structure is considered. Nevertheless, we can analyze the ensemble learning methods by using the replicated system with hierarchical structure induced by the nested quenched randomness.

There are two notable comments here. First, in the above discussion, we considered directly extrapolating the results for $n_1,n_2\in\sN$, obtained by considering only the trivial symmetry of the replicated system, to $n_1, n_2\to0$, which corresponds to the spirit of the replica symmetric (RS) solution. However, if the cost function (or the Hamiltonian) is non-convex, this extrapolation could lead to wrong results, and we have to include an additional structure of the replica symmetry breaking (RSB) by adding the replica symmetry breaking field \cite{parisi1989mechanism}, although we will not consider the effects of RSB in the following, since the loss functions we consider in this paper \eqref{eq: cost function SS} and \eqref{eq: cost function dKO} are both convex. Second, although we have focused on the limit $n_1,n_2\to0$ to evaluate the typical fluctuation over the algorithmic randomness $D^\ast$, one can consider the large deviation with respect to $D^\ast$ by analyzing the system with non-zero $n_2$. This can be understood by considering the generating function, similarly to the conventional formula for the free energy $\E[\log Z]=\lim_{n\to0}\partial_n \E[Z^n]$. From this generating function, we can understand that the replicated system \eqref{eq: replicated system} arises naturally in the analysis of randomized algorithms. See Appendix \ref{appendix: generating function approach} for more details.

\section{Mean-field Picture}
\label{section: mean-field picture}
In this section, we present the asymptotic properties of the estimators \eqref{eq: estimator SS} and \eqref{eq: cost function dKO} obtained by analyzing the replicated system \eqref{eq: replicated system}. In particular, we will present the results and implications obtained by considering the following generalized $K$-point correlation function defined as 
\begin{equation}
    \begin{split}
        \left\langle 
            \prod_{k=1}^K \psi(\theta_{\a^{(k)}, i})
        \right\rangle_\n
        \quad 
        a_1^{(1)}\neq a_1^{(2)}\neq\dots a_1^{(K)},
    \end{split}
    \label{eq: generalized K point correlation function}
\end{equation}
where $\psi$ is arbitrary as long as the following expectation is convergent. It is described by effective single-body problems whose behavior is determined by a small finite number of scalar quantities as the solution of nonlinear equations, which we refer to as self-consistent equations. Since the derivation is completely analogous to the 1RSB solution of the linear models in the limit where the breaking parameter goes to zero, we skip the derivation of them. Readers interested in their derivation can refer to the textbook \cite{engel2001statistical}.

\subsection{SS}
We first present the results for the estimator of SS \eqref{eq: estimator SS}. For this, let us define a single-body description of the estimator \eqref{eq: estimator SS} given as follows
\begin{equation}
    \begin{aligned}
        \wsfss(\xi,w_0,\eta) &= \argmin_{w\in\R}\frac{\hatqss}{2}w^2 - h_{\rm SS}w + \lambda|w|
        \\
        h_{\rm SS} &= \hatmss w_0 + \sqrt{\hatchiss}\xi + \sqrt{\hatvss}\eta
        \\
        &\xi,\eta\sim\gN(0,1), ~ w_0\sim p_{w_0},
    \end{aligned}
    \label{eq: effective single body description of SS}
\end{equation}
where $\hatqss, \hatmss, \hatchiss, \hatvss\in(0,\infty)$ are the parameters that are determined as a solution to the self-consistent equations defined in \eqref{eq: self-consistent ss f}-\eqref{eq: self-consistent ss v}, and $h_{\rm SS}$ is the effective local field to determine the random value of the estimator depending on $D,D^\ast$. 

Evaluating the generalized $K$-point correlation function \eqref{eq: generalized K point correlation function} at $n_1,n_2\to 0,~\beta\to\infty$, we obtain the following
\begin{equation}
    \E_{D}\left[
        \E_{D_{\rm SS}^\ast}\left[
            \psi(\hat{w}_i(D,D_{\rm SS}^\ast) )
        \right]^K
    \right]
    = \E_{\xi,w_0}\left[
        \E_{\eta}\left[
            \psi(\wsfss(\xi, w_0,\eta))
        \right]^K
    \right],
    \label{eq: average behavior SS}
\end{equation}
where the site-dependence disappears because all sites are treated equally in the current setup. This result indicates that the fluctuation with respect to the resampling conditioned by $D$ is effectively described by the Gaussian random variable $\eta\sim\gN(0,1)$, and the fluctuation with respect to the noise and the random explanatory variable is described by the Gaussian random variable $\xi\sim\gN(0,1)$. Schematically, this behavior can be represented as follows:
\begin{equation}
    h_{\rm SS} = 
    \remark{\substack{
        {\rm correlation~with}~\bm{w}_0
    }}{
        \hatmss w_0
    }
    + \remark{
        \substack{{\rm fluctuation} \\ {\rm from}~\{\bm{x}_\mu, \epsilon_\mu\}}
    }{
        \sqrt{\hatchiss} \xi
    }
    + \remark{
        \substack{{\rm fluctuation} \\ {\rm from}~D_{\rm SS}^\ast\mid D}
    }{
        \sqrt{\hatvss} \eta
    }.
\end{equation}
This implies that the average behavior of the estimator \eqref{eq: estimator SS} can be written as follows:
\begin{equation}
    \begin{split}
        \lim_{N\to\infty}\E_D\left[
            \phi\left(
                \E_{D_{\rm SS}^\ast}\left[
                    \psi(\hat{w}_{{\rm SS},i}(D,D_{\rm SS}^\ast))
                \right]
            \right)
        \right]
        \\
        =\E_{\xi,w_0}\left[
            \phi\left(
                \E_\eta\left[
                    \psi(\wsfss(\xi,w_0, \eta))
                \right]
            \right)
        \right],
    \end{split}
    \label{eq: average behavior SS general}
\end{equation}
where $\phi$ is arbitrary functions whose fluctuation can be determined by the moment method, and hence, TPR and FDR are written as follows
\begin{align}
    &\Pi_{{\rm SS}}(\xi, w_0) = \E_{\eta}\left[\1\left(\wsfss(\xi,w_0,\eta)\neq 0\right)\right],
    \label{eq: selection probability SS}
    \\
    &{\rm TPR} = \E_{\xi\sim\gN(0,1), w_0\sim\phi_+}\left[\1\left(
        \Pi_{{\rm SS}}(\xi, w_0) > \Pi_{\rm th}
    \right)\right],
    \label{eq: TPR SS}
    \\
    &{\rm FDR} = \frac{
        (1-\rho)\E_{\xi\sim\gN(0,1)}\left[\1\left(
            \Pi_{{\rm SS}}(\xi, 0) > \Pi_{\rm th}
        \right)\right]
    }{
        \E_{\xi\sim\gN(0,1), w_0\sim p_{w_0}}\left[
            \1\left(
                \Pi_{{\rm SS}}(\xi, w_0) > \Pi_{\rm th}
            \right)
        \right]
    }.
    \label{eq: FDR SS}
\end{align}

The parameters $\hatqss$, $\hatmss$, $\hatchiss$, $\hatvss$ are determined as the solution to the set of self-consistent equations. Let $\mu_{\rm B}=\lim_{N\to\infty}M_B/M$ be the resampling rate. Then, the self-consistent equations are given as follows:
\begin{align}
    f_1 &= \E_{c \sim{\rm Poi}(\mu_B)}\left[
        \frac{c}{1+\chi_{\rm SS} c}
    \right],
    ~
    f_2 = \E_{c \sim{\rm Poi}(\mu_B)}\left[
        \left(\frac{c}{1+\chi_{\rm SS} c}\right)^2
    \right],
    \label{eq: self-consistent ss f}
    \\
    \hatqss &= \hatmss = \alpha f_1,
    \label{eq: self-consistent ss qhat}
    \\
    \hatchiss &= \alpha f_1^2(q_{\rm SS}-2m_{\rm SS}+\rho+\Delta),
    \label{eq: self-consistent ss chihat}
    \\
    \hatvss &=\alpha\left(
        (f_2 - f_1^2)(q_{\rm SS}-2m_{\rm SS}+\rho+\Delta) + v_{\rm SS} f_2
    \right),
    \label{eq: self-consistent ss vhat}
    \\
    q_{\rm SS} &= \E_{\xi\sim\gN(0,1),w_0\sim p_{w_0}}\left[\E_{\eta\sim\gN(0,1)}[\wsfss]^2\right],
    \label{eq: self-consistent ss q}
    \\
    m_{\rm SS} &= \E_{\xi,\eta\sim\gN(0,1), w_0\sim p_{w_0}}\left[w_0 \wsfss\right],
    \label{eq: self-consistent ss m}
    \\
    \chi_{\rm SS} &= \E_{\xi,\eta\sim\gN(0,1), w_0\sim p_{w_0}}\left[\frac{{\rm d}}{{\rm d}h_{\rm SS}} \wsfss\right],
    \label{eq: self-consistent ss chi}
    \\
    v_{\rm SS} &= \E_{\xi\sim\gN(0,1), w_0\sim p_{w_0}}\left[
        \E_{\eta\sim\gN(0,1)}\left[\wsfss^2\right] -\E_{\eta\sim\gN(0,1)}\left[\wsfss\right]^2
    \right],
    \label{eq: self-consistent ss v}
\end{align}
where $c\sim {\rm Poi}(\mu_{\rm B})$ corresponds to the number of times that each data point appears in the bootstrap data $D_{\rm SS}^\ast$.

The parameters $(q_{\rm SS}, m_{\rm SS},\chi_{\rm SS},v_{\rm SS})$ that appear in the self-consistent equations are the order parameters which are related to the estimator \eqref{eq: estimator SS} as follows
\begin{align}
    q_{\rm SS} &= \lim_{N\to\infty}\frac{1}{N}\sum_{i=1}^N \E_{D_{\rm SS}^\ast}[\hat{w}_{{\rm SS},i}]^2,
    \\
    m_{\rm SS} &= \lim_{N\to\infty}\frac{1}{N}\sum_{i=1}^N \E_{D_{\rm SS}^\ast}[\hat{w}_{{\rm SS},i}] w_{0,i},
    \\
    \chi_{\rm SS} &=\lim_{N,\beta\to\infty}\frac{1}{N}\sum_{i=1}^N \E_{D_{\rm SS}^\ast}\left[
        \beta(\langle \hat{w}_{{\rm SS},i}^2\rangle - \langle \hat{w}_{{\rm SS},i}\rangle^2)
    \right],
    \\
    v_{\rm SS} &= \lim_{N\to\infty}\frac{1}{N}\sum_{i=1}^N \left(
        \E_{D_{\rm SS}^\ast}[\hat{w}_{{\rm SS},i}^2] - \E_{D_{\rm SS}^\ast}[\hat{w}_{{\rm SS},i}]^2
    \right).
\end{align}
The interpretation of these quantities is straightforward. That is, $q_{\rm SS}, m_{\rm SS}$ are the macroscopic second order statistics of averaged estimator $\E_{D^\ast}[\hat{\bm{w}}_{\rm SS}(D,D_{\rm SS}^\ast)]$ that characterizes the geometrical relation between  $\E_{D^\ast}[\hat{\bm{w}}_{\rm SS}(D,D_{\rm SS}^\ast)]$ and $\bm{w}_0$. Similarly, $v_{\rm SS}$ characterizes the variance of $\hat{\bm{w}}_{\rm SS}(D,D^\ast_{\rm SS})$ against $D^\ast_{\rm SS}$. Finally, $\chi_{\rm SS}$ represent the susceptibility of the estimator averaged over $D^\ast_{\rm SS}$

Therefore, by numerically or analytically analyzing the self-consistent equations, we can investigate the behavior of the estimator \eqref{eq: estimator SS}. Note that even though we are considering the large dimensional limit, the problem under analysis is now finite dimensional. This is a salient feature of the proportional asymptotic regime with a simple data model.

\subsection{dKO}
The behavior of the estimator \eqref{eq: estimator dKO} of dKO is also characterized by an effective single-body descriptions for $\hat{\bm{w}}_{\rm dKO}, \hat{\bwt}_{\rm dKO}$ given as follows:
\begin{align}
    \wsfko(\xi,w_0,\eta) &= \argmin_{w\in\R} \frac{\hatqko}{2}w^2 - h_{\rm dKO}w + \lambda|w|,
    \label{eq: effective estimator dKO w}
    \\
    \wtsfko(\tilde{\eta}) &= \argmin_{\wt\in\R} \frac{\hatqkko}{2}\wt^2 - \tilde{h}_{\rm dKO}\wt + \lambda|\wt|,
    \label{eq: effective estimator dKO wtilde}
    \\
    h_{\rm dKO} &= \hatmko w_0 + \sqrt{\hatchiko}\xi + \sqrt{\hatvko}\eta,
    \\
    \tilde{h}_{\rm dKO} &= \sqrt{\hatvkko}\tilde{\eta}.
\end{align}
Similarly to SS, the parameters $\hatqko, \hatqkko, \hatmko, \hatchiko$, $\hatvko, \hatvkko$ are determined as the solution to the set of the self-consistent equations \eqref{eq: self-consistent ko qhat}-\eqref{eq: self-consistent ko chik} that will be shown later. Moreover, the average behavior of the estimator \eqref{eq: estimator dKO} can be evaluated using them as follows:
\begin{equation}
    \begin{split}
        \lim_{N\to\infty}&\E_D\left[
            \phi\left(
                \E_{\tilde{X}}\left[
                    \psi(\hat{w}_i(D,\tilde{X}), \hat{\wt}_i(D,\tilde{X}))
                \right]
            \right)
        \right]
        \\
        & =\E_{\xi, w_0}\left[
            \phi\left(
                \E_{\tilde{\eta,\tilde{\eta}}}\left[
                    \psi(
                        \wsfko(\xi,w_0,\eta), \wtsfko(\tilde{\eta})
                    )
                \right]
            \right)
        \right],
     \end{split}
     \label{eq: average behavior dKO general}
\end{equation}
where $\phi, \psi$ are arbitrary functions as long as the above average is convergent. The site-dependence again disappears because all sites are treated equally in the current setup. The formula \eqref{eq: average behavior dKO general} indicates that the fluctuation with respect to the knockoff variable $\tilde{X}$ conditioned by $D$ is effectively described by the Gaussian random variables $\eta, \tilde{\eta}$, and the fluctuation with respect to the noise and the explanatory variable is described by the another Gaussian random variable $\xi$. All of this information is encoded in the effective local fields $h_{\rm dKO}, \tilde{h}_{\rm dKO}$. Schematically, this can be represented as follows:
\begin{align}
    h_{\rm dKO} &= 
    \remark{\substack{
        {\rm correlation~with}~\bm{w}_0
    }}{
        \hatmko w_0
    }
    + \remark{
        \substack{{\rm fluctuation} \\ {\rm from}~\{\bm{x}_\mu, \epsilon_\mu\}}
    }{
        \sqrt{\hatchiko} \xi
    }
    + \remark{
        \substack{{\rm fluctuation} \\ {\rm from}~\tilde{X} \mid D}
    }{
        \sqrt{\hatvko} \eta
    },
    \\
    \tilde{h}_{\rm dKO} &= 
    \remark{
        \substack{{\rm fluctuation} \\ {\rm from}~\tilde{X} \mid D}
    }{
        \sqrt{\hatvkko} \tilde{\eta}
    }.
\end{align}
It is noteworthy that the randomness coming from $D$ acts only on $\wsfko$ that describes the coefficient $\hat{\bm{w}}_{\rm dKO}$ for the true explanatory variable $\bm{x}_\mu$, while $\bwt_{\rm dKO}$ is determined solely by the knockoff variable $\tilde{X}$. This supports the intuition that the $\bwt_{\rm dKO}$ should be determined independently of the true data generation process.

Using the above result, TPR and FDR are written as follows:
\begin{align}
    &\begin{aligned}
        &\Pi_{\rm dKO}(\xi,w_0) 
        \\
        &= \E_{\eta,\tilde{\eta}}\left[
            \1\left(
                |\wsfko(\xi,w_0,\eta)| - |\wtsfko(\tilde{\eta})| > Z_{\rm th}
            \right)
        \right],
    \end{aligned}
    \label{eq: selection probability dKO}
    \\
    &{\rm TPR} = \E_{\xi\sim\gN(0,1), w_0\sim\phi_+}\left[\1\left(
        \Pi_{{\rm dKO}}(\xi, w_0) > \Pi_{\rm th}
    \right)\right],
    \label{eq: TPR dKO}
    \\
    &{\rm FDR} = \frac{
        (1-\rho)\E_{\xi\sim\gN(0,1)}\left[\1\left(
            \Pi_{{\rm dKO}}(\xi, 0) > \Pi_{\rm th}
        \right)\right]
    }{
        \E_{\xi\sim\gN(0,1), w_0\sim p_{w_0}}\left[
            \1\left(
                \Pi_{{\rm dKO}}(\xi, w_0) > \Pi_{\rm th}
            \right)
        \right]
    }.
    \label{eq: FDR dKO}
\end{align}

Finally, the self-consistent equations to determine $\hatqko$, $\hatqkko$, $\hatmko, \hatchiko$, $\hatvko, \hatvkko$ are given as follows
\begin{align}
    &\hatqko=\hatqkko=\frac{\alpha}{1+\chi_{\rm dKO}+\tilde{\chi}_{\rm dKO}},
    \label{eq: self-consistent ko qhat}
    \\
    &\hatchiko=\frac{\alpha}{(1+\chi_{\rm dKO}+\tilde{\chi}_{\rm dKO})^2}(q_{\rm dKO}-2m_{\rm dKO}+\rho+\Delta),
    \label{eq: self-consistent ko chihat}
    \\
    &\hatvko = \frac{\alpha}{(1+\chi_{\rm dKO}+\tilde{\chi}_{\rm dKO})^2}(v_{\rm dKO}+\tilde{v}_{\rm dKO}),
    \label{eq: self-consistent ko vhat}
    \\
    &\hatvkko = \hatchiko + \hatvko,
    \label{eq: self-consistent ko vkhat}
    \\
    &q_{\rm dKO} = \E_{\xi\sim\gN(0,1),w_0\sim p_{w_0}}\left[\E_{\eta\sim\gN(0,1)}[\wsfko]^2\right],
    \label{eq: self-consistent ko q}
    \\
    &m_{\rm dKO} = \E_{\xi,\eta\sim\gN(0,1),w_0\sim p_{w_0}}\left[\wsfko w_0\right],
    \label{eq: self-consistent ko m}
    \\
    &\chi_{\rm dKO} = \E_{\xi,\eta\sim\gN(0,1),w_0\sim p_{w_0}}\left[\frac{\rm d}{{\rm d}h_{\rm dKO}}\wsfko \right],
    \label{eq: self-consistent ko chi}
    \\
    &v_{\rm dKO} = \E_{\xi\sim\gN(0,1),w_0\sim p_{w_0}}\left[
        \E_{\eta\sim\gN(0,1)}[\wsfko^2]-\E_{\eta\sim\gN(0,1)}[\wsfko]^2
    \right],
    \label{eq: self-consistent ko v}
    \\
    &\tilde v_{\rm dKO} = \E_{\tilde{\eta}\sim\gN(0,1)}\left[\wtsfko^2\right],
    \label{eq: self-consistent ko vk}
    \\
    &\tilde \chi_{\rm dKO} = \E_{\tilde{\eta}\sim\gN(0,1)}\left[
        \frac{{\rm d}}{{\rm d}\tilde{h}_{\rm dKO}}\wtsfko
    \right].
    \label{eq: self-consistent ko chik}
\end{align}
Similarly to SS, the parameters $q_{\rm dKO}, m_{\rm dKO}, \chi_{\rm dKO}, v_{\rm dKO}, \tilde{v}_{\rm dKO}, \tilde{\chi}_{\rm dKO}$ are the order parameters for dKO which are related to the estimator \eqref{eq: estimator dKO} as follows
\begin{align}
    &q_{\rm dKO} = \lim_{N\to\infty}\frac{1}{N}\sum_{i=1}^N \E_{\tilde{X}}\left[\hat{w}_{{\rm dKO},i}\right]^2,
    \\
    &m_{\rm dKO} = \lim_{N\to\infty}\frac{1}{N}\sum_{i=1}^N \E_{\tilde{X}}\left[\hat{w}_{{\rm dKO},i}\right]w_{0,i},
    \\
    &\chi_{\rm dKO} = \lim_{N,\beta\to\infty}\frac{1}{N}\sum_{i=1}^N \E_{\tilde{X}}\left[
        \beta( \langle w_{{\rm dKO},i}^2\rangle - \langle w_{{\rm dKO},i}\rangle^2)
    \right],
    \\
    &v_{\rm dKO} = \lim_{N\to\infty}\frac{1}{N}\sum_{i=1}^N \left(\E_{\tilde{X}}\left[\hat{w}_{{\rm dKO},i}^2\right]-\E_{\tilde{X}}\left[\hat{w}_{{\rm dKO},i}\right]^2\right),
    \\
    &\tilde{v}_{\rm dKO} = \lim_{N\to\infty}\frac{1}{N}\sum_{i=1}^N \E_{\tilde{X}}\left[\hat{w}_{{\rm dKO},i}^2\right],
    \\
    &\tilde{\chi}_{\rm dKO} = \lim_{N,\beta\to\infty}\frac{1}{N}\sum_{i=1}^N \E_{\tilde{X}}\left[
        \beta( \langle \wt_{{\rm dKO},i}^2\rangle - \langle \wt_{{\rm dKO},i}\rangle^2)
    \right].
\end{align}
The interpretation of these quantities are similar to that of SS. Although $\tilde{v}_{\rm dKO}$ is seemingly characterized by only a second moment with respect to $\tilde{X}$, this is because its mean is zero $\E_{\tilde{X}}[\hat{\bwt}_{\rm dKO}]=\bm{0}$. Thus, $\tilde{v}_{\rm SS}$ can be interpreted as the variance of $\hat{\bwt}_{\rm dKO}$ with respect to $\tilde{X}$.

\subsection{Cross-checking with numerical experiments}

\begin{figure*}[t!]
    \centering
    \includegraphics[width=\textwidth]{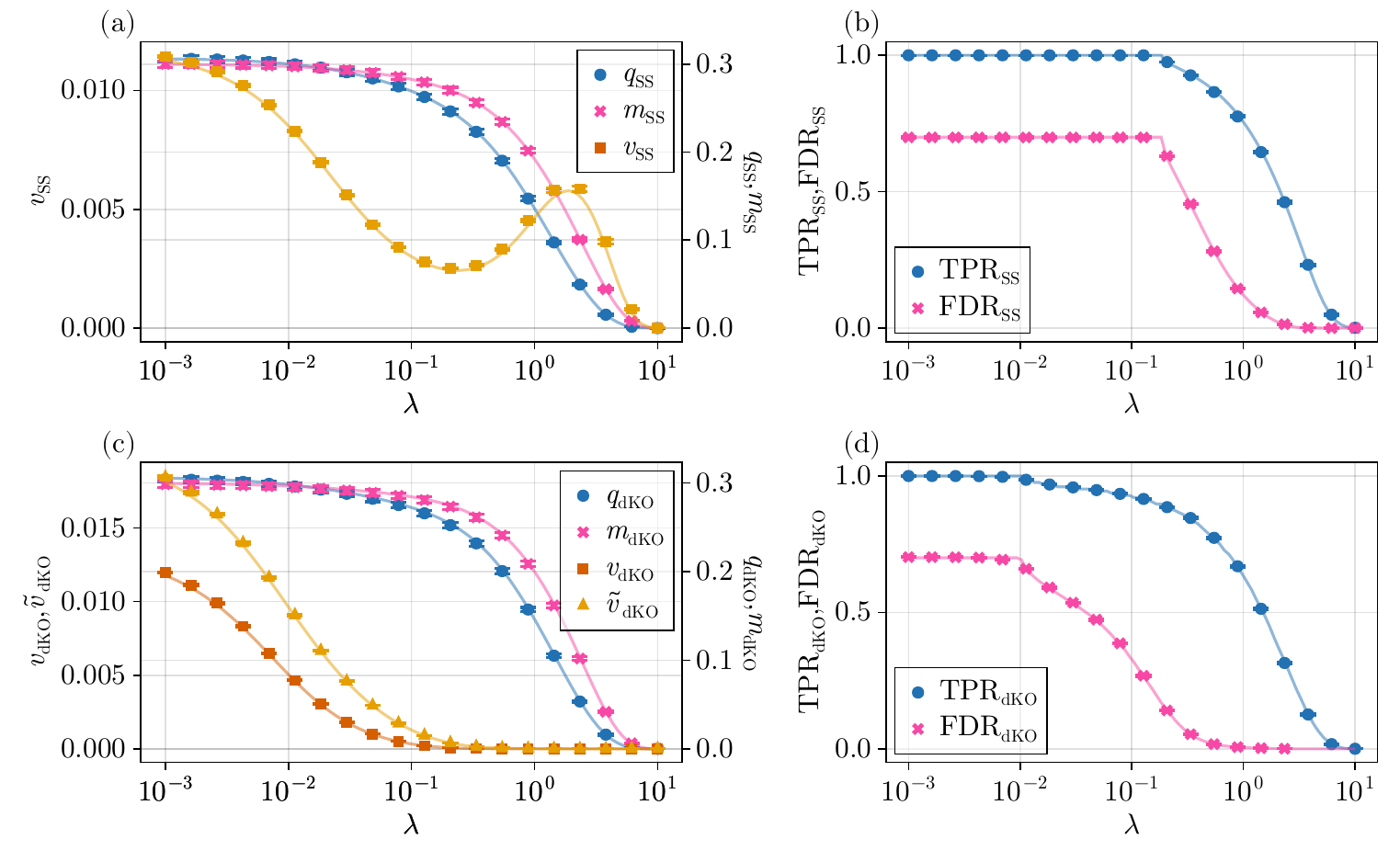}
    \caption{
        Comparison of theoretical predictions and experimental values for macroscopic quantities. (a)-(b) Results for SS. (c)-(d) Results for dKO. The markers with error bars represent the experimental values, where the error bars represent the standard error made by running experiments for several realizations of $D$. The solid lines represent the theoretical predictions. In all cases, parameters are set as $(\alpha,\rho,\Delta,\gamma_{\rm th}, \Pi_{{\rm th, dKO}} \Pi_{{\rm th,SS}})=(2.5,0.3,0.01,0.05,0.15,0.15)$.
    }
    \label{fig: comparison with experiments}
\end{figure*}

To check the validity of the formulas for the macroscopic quantities \eqref{eq: average behavior SS general} and \eqref{eq: average behavior dKO general} of SS and dKO, which are the most important results for evaluating the performance of these methods, we briefly compare them with numerical experiments of finite-size systems of $N=128$. In the experiments, the average over $D_{\rm SS}^\ast, \tilde{X}$ is approximated by the Monte Carlo method using $256$ realizations of them, and the estimators \eqref{eq: estimator SS} and \eqref{eq: estimator dKO} for each realization of $D, D_{\rm SS}^\ast, \tilde{X}$ are obtained by Glmnet \cite{friedman2010regularization}. The results are averaged over $512$ realizations of $D,D^\ast$ and the errors are evaluated as standard errors from them.

Fig. \ref{fig: comparison with experiments} shows the comparisons of the order parameters, TPR and FDR as functions of the regularization parameter $\lambda$. Since it is difficult to evaluate the second moment of the Boltzmann distribution, we omit the comparison for the susceptibilities $\chi_{\rm SS}$ and $\chi_{\rm dKO}$. Each panel corresponds to different algorithms or quantities. In all cases, they are in good agreement, demonstrating the validity of our analysis, despite $N$ being rather small.

\section{Performance Comparison}
\label{section: performance comparison}

In this section, we present the results for the comparison of the power, which is defined as the TPR at a given FDR, of each variable selection method. In addition to the standard SS with $\mu_B=M_B/M=1$ and dKO, we consider SS with $\mu_B=2$, which hopefully mitigate the decrease of unique data points in bootstrap data, and the vanilla knockoff which selects variables using single realization of $\tilde{X}$ as 
\begin{equation}
    \gS_{\rm KO}(D,\tilde{X}) = \{i \mid i\in[N], Z_i(D,\tilde{X})>Z_{\rm th}\}.
\end{equation}
The TPR and FDR for the vanilla knockoff can be easily evaluated using the single-body description \eqref{eq: effective estimator dKO w}-\eqref{eq: effective estimator dKO wtilde} as 
\begin{align}
    &{\rm TPR}_{\rm KO} = \E_{\xi\sim\gN(0,1),w_0\sim \phi_+}\left[
        \Pi_{\rm dKO}(\xi,w_0)
    \right],
    \\
    &{\rm FDR}_{\rm KO} = \frac{
        (1-\rho)\E_{\xi\sim\gN(0,1)}\left[
            \Pi_{\rm dKO}(\xi, 0)
        \right]
    }{
        \E_{\xi\sim\gN(0,1), w_0\sim p_{w_0}}\left[
            \Pi_{\rm dKO}(\xi, w_0)
        \right]
    }.
\end{align}
Recall that $\Pi_{\rm dKO}(\xi,w_0)$ is defined in \eqref{eq: selection probability dKO}.
As for the regularization parameter $\lambda$, we use $\lambda^\ast$ that minimizes the prediction error of the lasso estimators $\hat{\bm{w}}_{\rm SS}(D,D_{\rm SS}^\ast)$ or $\left(
    \hat{\bm{w}}_{\rm dKO}, \bwt_{\rm dKO}
\right)(D,\tilde{X})$.

\begin{figure}[t!]
    \centering
    \includegraphics[width=\linewidth]{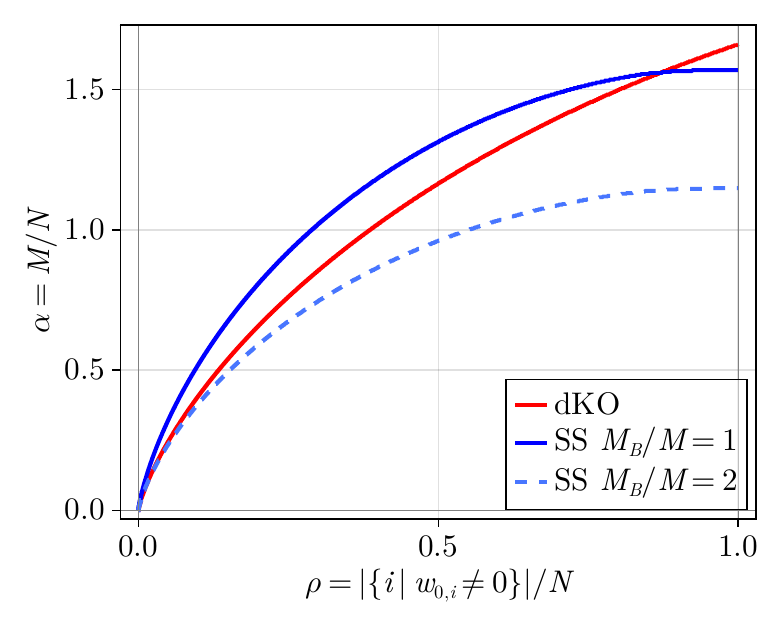}
    \caption{Comparison of perfect reconstruction limits \eqref{eq: perfect reconstruction limit} of each algorithm.}
    \label{fig: phase_boundary}
\end{figure}

Before directly comparing the power, we first consider the case when the measurement noise is absent, i.e., $\Delta=0$, which should be the easiest situation for a given $\alpha$. In this case, if the number of the unique data point in randomized data $D^\ast$ is sufficiently larger than some threshold, $\bm{w}_0$ should be perfectly recovered without an error in the sense that $\|\hat{\bm{w}}(D, D^\ast) -\bm{w}_0\|_2^2/N\to0$ \cite{kabashima2009typical,kabashima2010statistical,ganguli2010statistical}. Let $\tilde{\alpha}, \tilde{\rho}$ be the number of the unique data point in $D^\ast$ and the number of the non-zero elements of $\bm{w}_0$ normalized by the dimension of $\bm{\theta}$, respectively. Then, the condition for the perfect recovery is given as follows \cite{kabashima2009typical}
\begin{equation}
    \tilde{\alpha} > 2 (1-\tilde{\rho}) H(1/\sqrt{V}) - \tilde{\rho},
    \label{eq: perfect reconstruction limit}
\end{equation}
where $H(x) = \int_x^\infty \frac{e^{-\frac{z^2}{2}}}{\sqrt{2\pi}}dz$ and $V$ is the solution of the following equation
\begin{equation}
    \begin{split}
        \tilde{\alpha} V = 2(1-\tilde{\rho})&\left(
            (1+V)H\left(\frac{1}{\sqrt{V}}\right) - \sqrt{V}\frac{e^{-\frac{1}{2V}}}{\sqrt{2\pi}}
        \right) 
        \\
        &+ \tilde{\rho} (1+V),    
    \end{split}
\end{equation}

\begin{figure*}[t!]
    \centering
    \includegraphics[width=\linewidth]{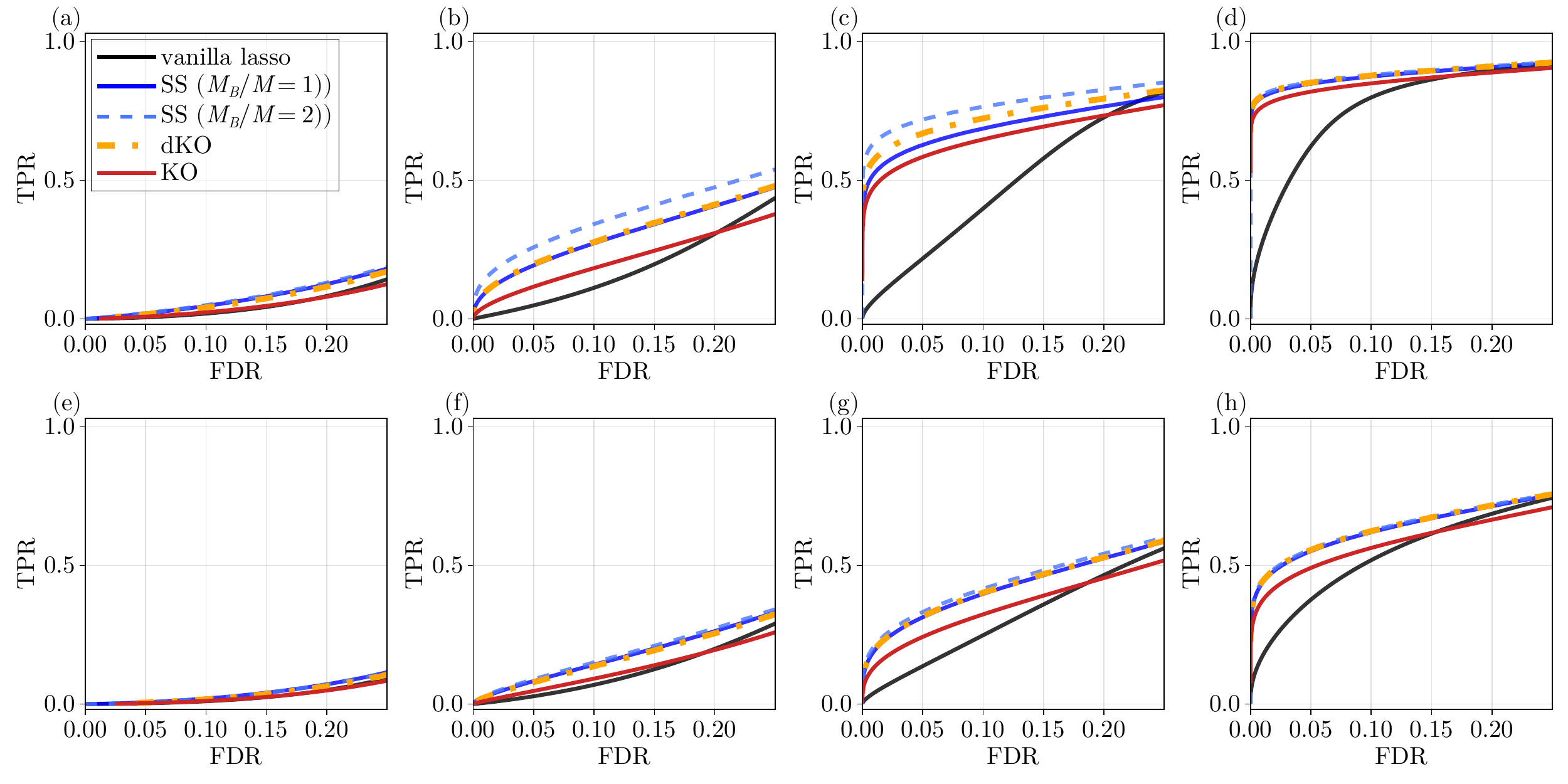}
    \caption{Comparisons of detection powers of several variable selection algorithms at several values of $(\alpha, \Delta)$. In all cases, $\rho$ is set as $0.5$. (a)-(d): Small noise case with $\Delta=0.01$. (e)-(h): Large noise case with $\Delta=0.1$. (a) and (e): $\alpha=0.36$. (b) and (f): $\alpha=0.63$. (c) and (g): $\alpha=1.12$. (d) and (h): $\alpha=2$.}
    \label{fig: power}
\end{figure*}

Although the above perfect reconstruction limit is defined in the noiseless case $\Delta=0$, it is expected that the estimators \eqref{eq: estimator SS} and \eqref{eq: estimator dKO} are highly correlated with $\bm{w}_0$ above the reconstruction limit and $\Delta$ is small. Therefore, the comparison of the reconstruction limit would provide insight into the performance of each algorithm. In SS with the resampling rate $\mu_B=M_B/M$ and dKO, $(\tilde{\alpha}, \tilde{\rho})$ is given as follows
\begin{equation}
    (\tilde{\alpha}, \tilde{\rho}) =\left\{
        \begin{array}{ll}
            \left((1-e^{-\mu_B})\alpha, \rho\right), & ~ {\rm SS},
            \\
            \left(\frac{\alpha}{2}, \frac{\rho}{2}\right), & {\rm dKO}.
        \end{array}
    \right.
\end{equation}
Using this relation, the comparison of the reconstruction limit is depicted in Fig. \ref{fig: phase_boundary}. It is clear that the reconstruction limit of dKO is located at a smaller $\alpha$ compared to the standard SS with $\mu_B=1$ as long as $\rho$ is not close to $1$, which is consistent with the argument made in the literature \cite{Ren2023derandomizing} that argues superiority of dKO over SS. However, with a resampling rate greater than $1$, the reconstruction limit shifts to an even smaller value of $\alpha$, suggesting the superiority of SS with $\mu_B>1$ over dKO.

Next, we will directly compare the detection power, which is defined as the TPR for a given FDR, of the variable selection algorithms. For SS, the detection power is plotted by changing $\Pi_{{\rm th, SS}}$, and for the vanilla KO and dKO, the detection power is plotted by changing $Z_{{\rm th}}$ while $\Pi_{{\rm th, dKO}}$ is fixed as $0.025$. Although the detection power of dKO depends both on $Z_{\rm th}$ and $\Pi_{{\rm th, dKO}}$, there was little difference between changing $Z_{\rm th}$ and changing $\Pi_{\rm th, dKO}$.  Also, the detection power of the vanilla lasso is also included for reference, where detection power is plotted as a function of the regularization parameter $\lambda$. In the vanilla lasso, the non-zero elements of the lasso estimator is directly used for estimating the support of $\bm{w}_0$. Finally, $\rho$ is set to $0.5$ in all cases.

In Fig. \ref{fig: power}, the comparison of the detection power is shown for several values of $(\alpha, \Delta)$. In all cases, the randomized algorithms including SS, the vanilla KO, and the dKO have higher detection power compared to the vanilla lasso, when FDR is close to zero, indicating the effectiveness of randomized algorithms. On the other hand, the superiority of the randomized algorithms depend on the conditions. In small noise cases (the upper panels (a)-(d)), it is clear that dKO provably outperforms both vanilla KO and the standard SS with $\mu_B=1$ at intermediate $\alpha$. However, the SS with $\mu_B=2$ is even higher detection power than dKO. This tendency is consistent with the behavior of reconstruction limit at $\Delta=0$. In large noise cases (the lower panels (e)-(h)), there are no difference between ensemble learning-based variable selection algorithms, although they are still better than the vanilla KO that depends on single realization of $\tilde{X}$ or the vanilla lasso. 

These results suggest that ensemble learning-based methods are superior to methods using a single realization of randomized data, and that in ensemble learning-based methods, it is recommended for each random dataset to contain as much information as possible.

\section{Summary and Discussion}
\label{section: summary and discussion}

In this paper, we presented the replica method for analyzing the ensemble learning-based variable selection method including SS and dKO, where the data size $M$ and the parameter dimension $N$ diverge at the same rate as $M,N\to\infty, M/N\to\alpha\in(0,\infty)$. The key to the analysis was the introduction of the replicated system \eqref{eq: replicated system} with a trivial hierarchical structure. Although the replicated system was defined heuristically in the main text, as shown in Appendix, it can be naturally derived by considering the generating function \eqref{eq: definition of the generating function} that describes the large deviation of the free energy with respect to the algorithmic randomness \eqref{eq: rate function}. 

It was shown that the behavior of the estimators for the randomized cost functions \eqref{eq: estimator SS} and \eqref{eq: estimator dKO} could be effectively described by the one-dimensional random variables \eqref{eq: effective single body description of SS} and \eqref{eq: effective estimator dKO w}-\eqref{eq: effective estimator dKO wtilde}, where the parameters are determined as a solution to the self-consistent equations \eqref{eq: self-consistent ss f}-\eqref{eq: self-consistent ss v} and \eqref{eq: self-consistent ko qhat}-\eqref{eq: self-consistent ko chik}. Consequently, the detection power formulas were also obtained as \eqref{eq: selection probability SS}-\eqref{eq: FDR SS} and \eqref{eq: selection probability dKO} and \eqref{eq: FDR dKO}.

By numerically evaluating the self-consistent equations, the comparison of the detection power was presented. It was shown that dKO outperforms the standard SS with resampling rate $\mu_B=1$ and the vanilla KO, while increasing the resampling rate in SS could further improve the detection power. 

The promising future directions include the extension of the analysis to more realistic data, because our analysis depends on iid assumption on the explanatory variables. Although our analysis suggests a superiority of the SS with over-sampled data, this could be an artifact of our setup because the performance of the lasso estimator is degraded by correlations between explanatory variables. Therefore, it is important to develop a framework to derive an effective single-body description that approximately describes the behavior in real data.

\begin{acknowledgment}

This study was supported by JSPS KAKENHI Grant No. 21K21310 and 23K16960, and Grant-in-Aid for Transformative Research Areas (A), “Foundation of Machine Learning Physics” (22H05117).
\end{acknowledgment}

\appendix

\section{Generating functional approach}
\label{appendix: generating function approach}
In this appendix, we briefly explain that the behavior of the replicated system \eqref{eq: replicated system} at $n_1\to0$ while keeping $n_2\neq0$ is related to the large deviation of the free energy \cite{ellis2007entropy} regarding the fluctuation of the algorithmic randomness. From this, it can also be understood that the replicated system, defined heuristically in the main text, is naturally derived from the analysis of the generating function, similarly to the conventional replica formula for the averaged free energy $\E[\log Z] = \lim_{n\to0}\partial_n\log \E[Z^n]$. As in Sec. \ref{section: replica method}, we use the notations $\bm{\theta}$ and $D^\ast$ for the parameters to be estimated and the algorithmic randomness, respectively, to treat SS and dKO in a unified manner.

For this, let us define the free energy density of the system 
\begin{equation}
    f^{(\beta)}(D,D^\ast) = \frac{-1}{N\beta}\log Z^{(\beta)}(D,D^\ast),
    \label{eq: free energy}
\end{equation}
and the probability density 
\begin{equation}
    p^{{(\beta)}}(f \mid D) = \E_{D^\ast}\left[
        \delta(f - f^{(\beta)}(D,D^\ast))
    \right],
    \label{eq: pdf of the free energy}
\end{equation}
that describes how the free energy fluctuates with respect to $D^\ast$ given $D$. Since the free energy is a macroscopic quantity, its fluctuation should be described by a rate function $s^{(\beta)}(f;D)$ as 
\begin{equation}
    p^{(\beta)} (f\mid D) \propto e^{Ns^{(\beta)}(f;D)}.
    \label{eq: rate function}
\end{equation}

For evaluating the rate function $s^{(\beta)}(f;D)$, it is convenient to define a generating function $g^{(\beta)}(\gamma;D), \gamma\in\R$ as 
\begin{equation}
    g^{(\beta)}(\gamma;D) \overset{\rm def}{=} \frac{1}{N}\log \int e^{-N\beta \gamma f}p^{(\beta)}(f\mid D) df.
    \label{eq: definition of the generating function}
\end{equation}
Once the generating function is evaluated for $\gamma\in\R$, one can obtain the rate function $s^{(\beta)}(f;D)$ from the Legendre transform of $g^{(\beta)}(\gamma;D)$ assuming the convexity of $g^{(\beta)}$. To see this, let us insert \eqref{eq: rate function} into \eqref{eq: definition of the generating function}, and evaluate the integral over $f$ at $N\gg1$ using the saddle point method as follows
\begin{align}
    g^{(\beta)}(\gamma;D) &= \frac{1}{N}\log \int e^{N(s^{(\beta)}(f;D) - \beta \gamma f)}df + o(N).
    \nonumber
    \\
    &\simeq \max_{f}[s^{(\beta)}(f;D)-\beta \gamma f].
\end{align}
Hence, the rate function is obtained as 
\begin{equation}
    s^{(\beta)}(f;D) = -\frac{\gamma^2}{\beta} \frac{\p}{\p\gamma}g^{(\beta)}(\gamma;D).
\end{equation}
In addition, $\lim_{\gamma\to0}g^{(\beta)}(\gamma;D)$ corresponds to the typical realization of $f^{(\beta)}(D,D^\ast)$ with respect to $D^\ast$ since $e^{-N\beta\gamma f}$ does nothing and the saddle point is $\max_f s^{(\beta)}(f;D)$.

In practice, the following expression of the generating function 
\begin{align}
    g^{(\beta)}(\gamma;D) &\overset{(a)}{=} \frac{1}{N}\log \E_{D^\ast}[(Z^{(\beta)}(D,D^\ast))^\gamma]
    \label{eq: generating function power of the partition function}
    \\
    &\overset{(b)}{=}\frac{1}{N}\E_D\left[
        \log \E_{D^\ast}[(Z^{(\beta)}(D,D^\ast))^\gamma]
    \right],
    \label{eq: generating function power of the partition function self averaging}
\end{align}
is convenient for the analysis. The equality (a) follows by simply inserting the definition of the density \eqref{eq: pdf of the free energy} into \eqref{eq: definition of the generating function}, and (b) follows by assuming the self-averaging of the generating function. Although the evaluation of the right-hand side of \eqref{eq: generating function power of the partition function self averaging} is still difficult since it contains the average of the log partition function, one can resort to the standard replica method as 
\begin{align}
        &\E_D\left[
            \log \E_{D^\ast}[(Z^{(\beta)}(D,D^\ast))^\gamma]
        \right] 
        = \lim_{n_1\to0}\frac{\partial}{\partial n_1}\log \tilde{\gZ}^{(\beta)}_{(n_1 \gamma)},
        \label{eq: replica formula for the generating function}
        \\
        & \tilde{\gZ}^{(\beta)}_{(n_1 \gamma)}=\E_D\left[
            \left(
                \E_{D^\ast}\left[
                    \left(Z^{(\beta)}(D,D^\ast)\right)^\gamma
                \right]
            \right)^{n_1}
        \right],
\end{align}
where the right-hand side of the above formula is evaluated by continuing the formula for $n_1\in\sN$ to $n_1\in\R$ following the standard prescription of the replica method. The important point is that, for $\gamma=n_2\in\sN$, $\tilde{\gZ}_{(n_1 \gamma)}^{(\beta)}$ is equal to $\gZ_\n^{(\beta)}$, which is the partition function of the replicated system \eqref{eq: replicated system}. Therefore, the replicated system \eqref{eq: replicated system} naturally appears from the analysis of the fluctuation over $D^\ast$ and the limit $n_2\to0$ corresponds to the analysis of typical fluctuation of $D^\ast$.


\bibliography{main}
\bibliographystyle{jpsj}

\end{document}

%% file: math_commands.tex
\usepackage{amsmath,amsfonts,bm}

\def\1{\bm{1}}

\DeclareMathAlphabet{\mathsfit}{\encodingdefault}{\sfdefault}{m}{sl}
\SetMathAlphabet{\mathsfit}{bold}{\encodingdefault}{\sfdefault}{bx}{n}

\def\gL{{\mathcal{L}}}

\def\gN{{\mathcal{N}}}

\def\gS{{\mathcal{S}}}

\def\gZ{{\mathcal{Z}}}

\def\sN{{\mathbb{N}}}

\newcommand{\E}{\mathbb{E}}

\newcommand{\R}{\mathbb{R}}

\DeclareMathOperator*{\argmin}{arg\,min}